\documentclass[11pt,letterpaper]{article}
\usepackage[T1]{fontenc}
\usepackage{lmodern}
\usepackage[margin=1in]{geometry}
\usepackage{amsmath,amssymb,amsthm,booktabs,enumitem,microtype}
\usepackage[hidelinks]{hyperref}
\usepackage{xurl}
\usepackage{titling,needspace}
\hypersetup{
 pdftitle={Vertex operator algebra realizations of untwisted Haagerup-Izumi centers},
 pdfauthor={David E. Evans and Andrew Riesen},
 pdfsubject={Research announcement, 23 October 2026}
}
\setlist[enumerate]{leftmargin=*,itemsep=4pt,topsep=4pt}
\newtheorem{theorem}{Theorem}
\newcommand{\HH}{\mathcal H}
\newcommand{\DD}{\mathcal D}
\newcommand{\PP}{\mathcal P}
\newcommand{\MM}{\mathcal M}
\newcommand{\one}{\mathbf1}
\newcommand{\Z}{\mathbb Z}
\newcommand{\C}{\mathbb C}
\newcommand{\e}{\mathrm e}
\DeclareMathOperator{\Rep}{Rep}

\DeclareMathOperator{\FPdim}{FPdim}
\DeclareMathOperator{\Ext}{Ext}

\DeclareMathOperator{\Vect}{Vec}
\DeclareMathOperator{\MD}{MD}
\title{Vertex operator algebra realizations of\\
untwisted Haagerup--Izumi centers}
\author{David E. Evans \and Andrew Riesen}
\date{October 23, 2026}

\begin{document}
\maketitle

\begin{abstract}
For every odd integer $N\geq3$, we construct an untwisted
Haagerup--Izumi fusion category $\mathcal H_N$ on the negative spherical
branch and a simple, self-contragredient, regular, $C_2$-cofinite vertex
operator algebra $E_N$ of central charge $-16N$, together with an equivalence of modular fusion categories:
\[
 \operatorname{Rep}(E_N)\simeq Z(\mathcal H_N).
\]

The VOA $E_N$ is rational, has exactly $N^2+3$ irreducible modules,
and is lower bounded and integer graded, with least conformal degree
$-(N-1)/2$.
The proof begins with a CFT-type cyclic Virasoro extension of central
charge $-5N$, whose regularity we establish. Two categorical
condensations recover $\mathcal H_N$, and a finite extension by
shifted lattice modules removes the auxiliary pointed factor.
We identify the complete ribbon modular data of $Z(\mathcal H_N)$ with
$\mathrm{MD}_{-}(\mathbb Z/N;\mathbb Z/(N^2+4),b\mapsto b^2/(N^2+4))$
in the notation of Gannon--Schopieray--Yadav, for every odd $N$.
\end{abstract}

\section{Introduction and main result}
\label{sec:intro}

The point of this article is to announce that the authors, with the
assistance of AI, have found an infinite family of VOAs realizing centers of untwisted cyclic
Haagerup--Izumi categories. The construction applies to every odd
order and lies on the negative spherical branch. More precisely,
we announce the following result.

\begin{theorem}[{\cite[Theorems 1.2, 30.2 and 31.1]{ER}}]
\label{thm:main}
For every odd integer $N\geq3$, put $G=\Z/N$ and $\mu=N^2+4$.
There are an untwisted cyclic Haagerup--Izumi spherical fusion category
$\HH_N$, whose Drinfeld center $Z(\HH_N)$ is a modular tensor category,
and simple, self-contragredient, regular, $C_2$-cofinite
vertex operator algebras $E_N$ and $W_N$ with:
\[
\begin{aligned}
 c(E_N)&=-16N,&
 \Rep(E_N)&\simeq_{\mathrm{ribbon}}Z(\HH_N),\\
 c(W_N)&=-5N,&
 \Rep(W_N)&\simeq_{\mathrm{ribbon}}\PP_N\boxtimes Z(\HH_N),
\end{aligned}
\]
where $\PP_N$ is a nondegenerate pointed ribbon category of rank $2^N$.
Both VOAs are rational. The VOA $E_N$ has exactly $N^2+3$
irreducible modules and least conformal degree $-(N-1)/2$;
it is lower bounded and integer graded, but not of CFT type.
The VOA $W_N$ is of CFT type and $(W_N)_1=0$.

The category $\HH_N$ has simple objects $g$ and $g\rho$, $g\in G$,
with spherical dimensions $\dim g=1$ and $\dim(g\rho)=-t$, and
$\FPdim(g\rho)=\delta$, where:
\[
 t=\frac{\sqrt\mu-N}{2},\qquad
 \delta=t^{-1}=\frac{N+\sqrt\mu}{2}.
\]
With $S_{\one,\one}=1/\dim\HH_N$ and $T^{\mathrm{rib}}$ the
ribbon twist matrix, the complete modular data of $Z(\HH_N)$ are:
\[
 (S,T^{\mathrm{rib}})
 =\MD_-\bigl(\Z/N;\Z/\mu,\ b\longmapsto b^2/\mu\bigr)
\]
in the conventions of Gannon--Schopieray--Yadav
\cite[Appendix F.1]{GSY}, including composite and nonsquarefree $\mu$.
\end{theorem}

Here $\HH_N$ is the category defined by the two condensation
algebras below. The ribbon equivalence is obtained from these
algebras; the modular-data comparison is a separate calculation.
Equivalence or categorical Morita equivalence with an independently
prescribed Haagerup--Izumi category remains a further question.

We outline the proof and its principal ideas below. Before doing so,
we describe the role of AI in the article and in the mathematical
investigation.

\subsection*{AI disclosure}
This announcement and the accompanying full manuscript were
prepared with extensive use of generative AI, including drafting
and editing. AI produced substantial portions of the initial
arguments and exposition. The authors directed the work and
revised the arguments and text; they take responsibility for the
mathematical claims and final version. Independent human
verification of the complete argument remains pending.

\subsection*{AI methodology}
The main AI tools used have been GPT Pro 6, Astra 6, Sol 6.1, Sol 5.6 and Opus 5.5.
Here is the order of events leading to the result. For the sake of
brevity, we omit the many prompts required to reorient the AI
properly.
\begin{enumerate}
\item \textbf{Initial discussions and the search for a uniform construction.}
Following the second author's announcement of
\href{https://arxiv.org/abs/2610.13933v1}{the Haagerup VOA at central charge eight}~\cite{Riesen},
the authors brought together that work and the first author's
\href{https://arxiv.org/abs/2610.18512v1}{investigations of S-folds}~\cite{SFolds} and began collaborating.
Our aim became a uniform construction for an infinite Haagerup--Izumi family.
AI-assisted exploration of possible constructions was guided by character
constraints, relevant literature, and other mathematical consistency tests.

\item \textbf{The twisted family.}
The first author's
\href{https://arxiv.org/abs/2610.26341v1}{subsequent VOA construction}~\cite{Evans},
together with his conjecture concerning its extension to an infinite family,
supplied a concrete starting point. Using this manuscript as input, both authors
directed ChatGPT to investigate the conjecture and develop arguments uniformly
in the parameter $N$. These arguments provided the basis for the further
development of the construction and our joint work on the manuscripts.

\item \textbf{Isolating the twisted doubles.}
A further idea was to adapt a reconstruction argument from
\href{https://arxiv.org/abs/2410.00809}{\emph{Orbifolds of Pointed Vertex Operator Algebras I}}~\cite{GR},
joint work of Gannon and the second author, to remove the additional pointed
factor and isolate the desired twisted doubles. We asked ChatGPT to develop
this approach, leading to the construction presented in our joint manuscript,
\href{https://hexagonmath.org/2610.00490v1}{\emph{Vertex Operator Algebra Realizations of Twisted Haagerup--Izumi Doubles}}~\cite{Twisted}.

\item \textbf{The untwisted family.}
Building on the twisted construction, we then directed the investigation toward
the untwisted case. Further mathematical guidance and exchanges with ChatGPT led
to the construction announced here, which we have continued to develop and
revise jointly.

\item \textbf{Review and announcement.}
AI systems, including ChatGPT and Claude, played a substantial role in
developing, examining, and revising the arguments, as well as in drafting the
manuscripts. Throughout the investigation, we supplied mathematical ideas,
references, constraints, and corrections, and used the resulting arguments to
guide subsequent steps. Both authors have been checking the arguments and
revising the exposition. The scope of the constructions and our increasing
confidence in them have motivated this announcement, through which we hope to
bring these developments to the attention of the mathematical community while
verification continues. We take responsibility for the mathematical claims and
welcome scrutiny of the arguments.
\end{enumerate}
\paragraph{Working manuscripts.}
The current working manuscript~\cite{ER} and the companion manuscript
on twisted Haagerup--Izumi doubles~\cite{Twisted}, together with their
LaTeX sources, are available on GitHub at:
\[
\text{\url{https://github.com/Doctorfm98/Haagerup_Izumi_VOA}}.
\]
Both manuscripts are also available on Hexagon.
As described in the methodology, the work remains ongoing.
The authors will continue refining the arguments and exposition
to make the proofs easier for readers to follow and verify.
Corrections and suggestions are welcome.
\subsection*{Background and method}
The centers of Haagerup--Izumi fusion categories have long been a
test case for the realization of modular tensor categories by
vertex operator algebras. Izumi's endomorphism and tube-algebra
constructions \cite{Izumi}, followed by the modular-data and
nonunitary constructions of the first author and Gannon \cite{EG11,EG17},
made this problem explicit. Recent work of Huang \cite{Huang} and
Gannon--Schopieray--Yadav \cite{GSY} establishes the existence of
cyclic Haagerup--Izumi categories at every odd order. We construct
regular vertex operator algebras realizing centers of such
categories uniformly in the odd order, with the negative
spherical structure.

The spherical structure is part of the result. Recent $c=8$
announcements concern the positive branch at $N=3$
\cite{Riesen,HuangThree} and $N=5$ \cite{HuangFive}.
The family constructed here is nonunitary: its noninvertible
Haagerup--Izumi objects have the negative root of
$d^2=1+Nd$ as their spherical dimension. Their
Frobenius--Perron dimension is the positive root. This distinction
persists in the ribbon structure of the center.

The starting point is the cyclic Virasoro sum proposed by the first author
\cite[Section 6.2]{Evans} in the S-fold setting, following the
duality-wall gluing of Creutzig--Gaiotto \cite{CG}. We give this
cyclic sum its VOA structure, prove regularity, and determine all
of its simple modules. The lattice realization also supplies the
algebra objects from which the Haagerup--Izumi category and its
ribbon center are recovered. A finite extension by shifted $E_7$
lattice modules then realizes the center itself. Detailed proofs
appear in the accompanying full manuscript \cite{ER}.

\Needspace{10\baselineskip}
\paragraph{Conventions.}
All categories are over $\C$. Dimensions are spherical dimensions
unless marked $\FPdim$, and $Z(\HH)$ carries the ribbon structure
induced by the specified spherical structure on $\HH$.
The notation $G$, $\mu$, $t$ and $\delta$ remains as in
Theorem~\ref{thm:main}; in particular, $t(t+N)=1$.
An \emph{untwisted cyclic Haagerup--Izumi category} has fusion rules:
\[
 gh=g+h,\qquad \rho g=(-g)\rho,\qquad
 \rho^2=\one\oplus\bigoplus_{g\in G}g\rho,
\]
and pointed subcategory tensor equivalent to $\Vect_G$ with
trivial associator. The remaining associativity data are part of
the construction.

We allow a VOA to be lower bounded and integer graded, with
finite-dimensional homogeneous spaces. \emph{Regularity} means
that every weak module is a direct sum of irreducible ordinary
modules. For the regular VOAs here, $\Rep(V)$ is the category of
grading-restricted generalized modules, consisting of finite
direct sums of ordinary simples.
CFT type means $V_n=0$ for $n<0$ and $V_0=\C\one$;
rationality means complete reducibility of admissible modules.
Self-contragredience means $V\simeq V^\vee$ for the stated conformal vector.
The condition of $C_2$-cofiniteness is:
\[
 C_2(V)=\operatorname{span}_{\C}\{a_{-2}b:a,b\in V\},
 \qquad \dim V/C_2(V)<\infty.
\]

\section{The cyclic VOA and its regularity}
\label{sec:construction}

For $0\leq i<N$, set:
\[
 \psi_i=\frac{t+i+1}{t+i},\qquad
 c_i=13-6(\psi_i+\psi_i^{-1}),\qquad
 h_i(r,s)=\frac{(\psi_i r-s)^2-(\psi_i-1)^2}{4\psi_i}.
\]
Let $L_i(r,s)$ denote the simple Virasoro module of central charge
$c_i$ and lowest weight $h_i(r,s)$. Indices are cyclic. The underlying
frame decomposition is:
\begin{equation}
 W_N=\bigoplus_{r_0,\ldots,r_{N-1}\in2\Z_{\geq0}+1}
             \ \bigotimes_{i=0}^{N-1}L_i(r_i,r_{i+1}).
 \label{eq:cyclic}
\end{equation}
The relations $\psi_i+\psi_{i-1}^{-1}=k_i$, where
$(k_0,k_1,\ldots,k_{N-1})=(N+2,2,\ldots,2)$, pair the relevant
generic Virasoro tensor categories by braid-reversed equivalences.
The generic tensor-category results \cite{CJORY,Koshida} and diagonal
gluing \cite{Gluing} then give a multiplication on
\eqref{eq:cyclic}, with nonzero allowed fusion channels.
Writing $r_i=2m_i+1$ gives:
\[
 h(\mathbf m)=\sum_i(m_i-m_{i+1})^2+Nm_0(m_0+1).
\]
This gives the CFT grading, finite-dimensional homogeneous spaces,
and $(W_N)_1=0$. A nonzero ideal contains a frame summand;
multiplication with its dual reaches the vacuum, proving
simplicity. The invariant-form criterion then gives
self-contragredience. Summing the $c_i$ yields $c(W_N)=-5N$.

\Needspace{13\baselineskip}
The lattice realization makes the regularity argument accessible.
Define the positive form $K_+$ by:
\[
 x^{\mathsf t}K_+x=\sum_i(x_i-x_{i+1})^2+Nx_0^2,
 \qquad x_N=x_0,
\]
and let $L=(\Z^N,2K_+)$, with basis $\beta_i$ and
$B=\sum_i\beta_i$. An explicit Virasoro frame in a conformal
deformation of $V_L$ realizes \eqref{eq:cyclic} as a conformal
subalgebra. The ambient conformal vector has the form:
\[
 \omega_L+\tfrac12B(-2)\one+\lambda D^2e^B,
 \qquad \lambda\ne0,
\]
where $D$ is translation and the normalization is fixed by the frame.
Recognition of the algebra multiplication, not just its frame
decomposition, is part of the embedding argument
\cite[Sections 4--10]{ER}.

The proof of regularity precedes the use of modular transformations.
First, strong generation and the Virasoro null relations make
$W_N/C_2(W_N)$ finite over the algebra of stress-tensor classes.
At any closed point of this commutative algebra, those same relations
transport the stress values around the cycle with multiplier $t^4$.
Since $0<t<1$, all stress values vanish. Finite generation then gives
nilpotence of the stress classes and $C_2$-cofiniteness
\cite[Section 3]{ER}.

Next, one proves:
\begin{equation}
 \Ext^1(W_N,S)=0
 \quad\text{for every simple grading-restricted generalized module }S.
 \label{eq:ext}
\end{equation}
Frame null relations reduce a possible nonsplit vacuum extension to
lowest weights zero and one. In weight zero, filtered Zhu relations
and two distinct ground characters give a semisimple two-dimensional
quotient; a separate mode argument lifts this splitting to the
extension. In weight one, the relations for proper cyclic intervals
and the closing generator force the defect to vanish, while nonzero
vacuum coupling excludes the remaining possibility. The argument
includes generalized lowest spaces and the $N=3$ endpoint cases
\cite[Sections 12--16 and Appendices A--H]{ER}.

Equation~\eqref{eq:ext} makes the vacuum projective. In the
Grothendieck--Verdier tensor category supplied by VOA duality
\cite{ALSW}, the evaluation onto the simple vacuum therefore splits.
Finiteness of the abelian category and right exactness of tensor
product give exponential bounds on tensor-power lengths, and hence
the moderate-growth bounds required by the rigidity criterion of
Etingof--Penneys \cite[Theorem 1.1]{EP}.
Rigidity and projectivity of the unit make every object projective,
so the category is semisimple. The passage from cyclic weak modules
to arbitrary weak modules proves regularity
\cite[Section 17]{ER}. The resulting CFT-type VOA $W_N$ is strongly
rational in the usual sense, and its representation category is
modular \cite{McRae}.

\section{Simple modules and reconstruction of the center}
\label{sec:center}

The parity currents $J_\epsilon$, indexed by
$P=(\Z/2)^N$, are obtained by allowing
$r_i\equiv1+\epsilon_i\pmod2$ in the cyclic frame sum.
Their monodromy pairing is
$(-1)^{\epsilon^{\mathsf t}K_+\eta}$.
Since $\det K_+=N^2$ is odd, they form a nondegenerate pointed
subcategory $\PP_N$. Consequently:
\[
 \Rep(W_N)\simeq\PP_N\boxtimes\DD_N,
 \qquad \DD_N=\PP_N',
\]
where the prime denotes the braided centralizer.

Two lattice forms control the other simple modules. The periodic
form $K_+$ has determinant $N^2$. Replacing its closing difference
by $(x_{N-1}+x_0)^2$ gives an antiperiodic form $K_-$ with
determinant $\mu=N^2+4$. The corresponding rational frame grids
yield the periodic and exceptional families. Simplicity,
multiplicity one, and uniqueness on these grids are established
before the modular comparison. Refined traces retain the separate
frame variables, as formal germs along their diagonal, so that
exceptional labels remain distinguishable when $\mu$ is composite
\cite[Sections 19--25]{ER}. Scalar conformal-weight residues alone
would not suffice.

\Needspace{13\baselineskip}
The complete list in $\DD_N$ is as follows. Write
$\e(x)=\exp(2\pi i x)$; bracketed labels denote orbits under sign.
\begin{center}
\renewcommand{\arraystretch}{1.16}
\begin{tabular}{@{}lccc@{}}
\toprule
Simple objects & Number & Spherical dimension & Ribbon twist\\
\midrule
$\one$ & $1$ & $1$ & $1$\\
$T$ & $1$ & $t^2$ & $1$\\
$Y_{[x,y]}$ & $(N^2-1)/2$ & $1+t^2$ & $\e(xy/N)$\\
$X_{[a]}$ & $(\mu-1)/2$ & $-Nt$ & $\e(Na^2/\mu)$\\
\bottomrule
\end{tabular}
\end{center}
Here $(x,y)$ runs through the nonzero elements of $(\Z/N)^2$
modulo sign, and $a$ through $(\Z/\mu)\setminus\{0\}$ modulo sign.
Thus $\DD_N$ has $N^2+3$ simple objects, whereas $\HH_N$ has $2N$.

Exhaustion of the module list uses more than candidate characters.
The lattice extension gives a connected \'etale algebra
$A=\one\oplus T$ in $\DD_N$; its categorical index is shown to be
nonzero before invoking separability. Branching gives $N^2$
invertible $A$-modules and one further simple module $R$.
Write $F_A(T)=\one\oplus kR$ for free-module induction, with
$k$ a nonnegative integer, and let $n$ be the common positive
induction multiplicity of the known exceptional modules.
Dimension and modular-coefficient identities give:
\[
 n^2(N^2+4)=N^2k^2+4,\qquad
 \sum_{\text{additional }Z}n_Z^2=\frac{2(1-k^2)}{N^2+4}.
\]
Nonnegativity forces $k\leq1$, and the first equation excludes
$k=0$. Hence $k=n=1$ and there are no additional simples
\cite[Sections 26--28]{ER}.

The resulting category $\MM_N=(\DD_N)_A$ has near-group rules:
\[
 eR=Re=R,\qquad
 R^2=\bigoplus_{e\in(\Z/N)^2}e\ \oplus N^2R.
\]
Its local-module category is the hyperbolic pointed ribbon
category on $G\times\widehat G$, with quadratic form
$q(g,\chi)=\chi(g)$. The central embedding of the reverse local
category into $Z(\MM_N)$, whose underlying functor is fully
faithful, supplies the central algebra for a second condensation.
This retains the half-braiding and multiplication needed to
recover the HI category.

Choose one cyclic Lagrangian factor $H$ of
$G\times\widehat G$ and its regular algebra. Condensation by this
central algebra identifies the invertibles modulo $H$.
The induced $R$ has endomorphism algebra a twisted group algebra
of the cyclic group $H$. Since $H^2(H,\C^\times)=0$, it splits
into $N$ distinct simple summands. The deck action and Frobenius
reciprocity give exactly the HI fusion rules stated in
Section~\ref{sec:intro}. The complementary Lagrangian factor
supplies a tensor equivalence from $\Vect_G$ to the pointed
subcategory, proving that its associator is trivial
\cite[Section 29]{ER}.

If $A'$ denotes the composite algebra in $\DD_N$, the construction
therefore defines:
\[
 \HH_N=(\DD_N)_{A'},\qquad
 (\DD_N)_{A'}^0\simeq\Vect.
\]
The composite algebra $A'$ is connected \'etale, has trivial
twist and nonzero spherical dimension, and is Lagrangian.
The center theorem \cite[Corollary 4.1]{DMNO} therefore gives
$\DD_N\simeq Z(\HH_N)$ as braided categories.
The induced spherical structure satisfies
$\dim_{\HH_N}M=\dim_{\DD_N}M/\dim_{\DD_N}A'$
\cite[Theorems 1.17--1.18 and Remark 1.19]{KO}; in particular
$\dim\rho=-t$. The center equivalence preserves these dimensions,
which determine the pivotal structure on simple objects, and
hence preserves the twists. This gives the ribbon equivalence in
Theorem~\ref{thm:main}.
Since $\HH_N$ is a spherical fusion category over $\C$, its
Drinfeld center $Z(\HH_N)$ is modular by the general center theorem
\cite[Section 2.9]{DGNO}.

\section{An explicit VOA for the center}
\label{sec:extension}

The pointed factor can be removed at the VOA level. Define the
invertible linear map $\sigma:P\to P$ by:
\[
 (\sigma\epsilon)_i=\epsilon_i+\epsilon_{i+1}+\epsilon_0
 \quad\text{in }\mathbb F_2.
\]
In these coordinates, the self-braiding, twist and dimension of
$J_\epsilon$ are respectively
$i^{|\sigma\epsilon|}$, $(-i)^{|\sigma\epsilon|}$ and
$(-1)^{|\sigma\epsilon|}$.

Let $K_7$ be the $E_7$ root lattice and choose its minuscule vector
$u\in K_7^*\setminus K_7$, with $u^2=3/2$. Set:
\[
 U=(V_{K_7},\omega_{K_7}+u(-2)\one),\qquad M=V_{K_7+u}.
\]
Conformal shifting \cite{DM} gives $c(U)=-11$.
The shifted $M$ has lowest weight $-3/4$, self-braiding $-i$,
twist $i$, and spherical dimension $-1$. Thus $\Rep(U)$
is the reverse ribbon category of one parity factor.
With $M^0=U$ and $M^1=M$, the required extension is:
\begin{equation}
 E_N=\bigoplus_{\epsilon\in P}
       J_\epsilon\otimes\bigotimes_{i=0}^{N-1}M^{(\sigma\epsilon)_i}.
 \label{eq:extension}
\end{equation}
Paired self-braidings and twists cancel. The diagonal algebra
is connected, commutative and separable, and the extension
correspondence \cite{CMY,CKM} gives the VOA structure.
Its local modules are precisely $\DD_N$:
\[
 c(E_N)=-5N-11N=-16N,\qquad
 \Rep(E_N)\simeq_{\mathrm{ribbon}}\DD_N.
\]

The grading and finiteness assertions require separate arguments.
If $T(\epsilon)$ counts transitions in the cyclic binary word
$\epsilon$, the minimum degree of its summand in
\eqref{eq:extension} is:
\[
 h_{\min}(E_{N,\epsilon})=
 \begin{cases}
 -T(\epsilon)/2,&\epsilon_0=0,\\
 T(\epsilon),&\epsilon_0=1.
 \end{cases}
\]
On an odd cycle, the maximum transition count is $N-1$, attained
with $\epsilon_0=0$. All higher degrees differ by integers.
This gives the precise least degree $-(N-1)/2$.

To prove $C_2$-cofiniteness, use the same underlying vertex
algebra as the conformal base, but restore its unshifted lattice
conformal vectors:
\[
 V^{\mathrm{un}}=W_N\otimes V_{K_7}^{\otimes N}.
\]
This algebra is of CFT type and $C_2$-cofinite.
Equation~\eqref{eq:extension} is a finite sum of its irreducible
weak modules. Each is $C_2$-cofinite as a module by
\cite[Proposition 5.2]{ABD}. Set:
\[
 C_2^{V^{\mathrm{un}}}(E_N)
 =\operatorname{span}_{\C}\{v_{-2}e:v\in V^{\mathrm{un}},\ e\in E_N\}.
\]
This space has finite codimension, and:
\[
 C_2^{V^{\mathrm{un}}}(E_N)\subseteq C_2(E_N),\qquad
 \dim E_N/C_2(E_N)<\infty.
\]
The conformal vectors need not agree for this inclusion:
the $C_2$-spaces depend only on the vertex operators.

\Needspace{13\baselineskip}
For regularity, the actual conformal base
$V=W_N\otimes U^{\otimes N}$ is nonnegatively graded and
$C_2$-cofinite, with semisimple module category.
Every cyclic weak $V$-module is grading-restricted generalized
\cite[Proposition 2.11]{McRae}; semisimplicity therefore makes
every weak $V$-module a direct sum of ordinary simples.
The extension correspondence on direct limits
\cite[Theorems 7.5 and 7.7]{CMY} identifies weak $E_N$-modules
with local algebra modules in this completion.
For such a module $Q$, the image of
$E_N\boxtimes_V X\to Q$, with $X\subset Q$ a finite
$V$-subobject, is an $E_N$-submodule finite as a $V$-object.
These images cover $Q$.
Separability makes each image semisimple, so $Q$ is a sum, and
hence a direct sum, of simple submodules. Every simple summand
is finite as a $V$-object and is ordinary for the common conformal
vector. Thus $E_N$ is regular for arbitrary weak modules.

Finally, $2u\in K_7$ makes $U$, and hence $V$, self-contragredient.
The $V$-summand of $E_N$ gives a $V$-summand of its
contragredient. Induction--restriction adjunction then gives a
nonzero map $E_N\to E_N^\vee$; simplicity makes it an isomorphism
\cite[Theorem 31.1]{ER}. These arguments establish regularity,
rationality, $C_2$-cofiniteness and self-contragredience for the
stated conformal vector of $E_N$.

\section{The modular data and the metric group}
\label{sec:modular}

The exceptional twists initially give the metric group
$(\Z/\mu,q_N)$, where $q_N(a)=Na^2/\mu$ in $\mathbb Q/\mathbb Z$.
There is a uniform isometry to the cyclic form
$q_0(b)=b^2/\mu$. Indeed, set:
\[
 u_N=\frac{N^2-N+2}{2}.
\]
The congruence $q_0(u_Na)=q_N(a)$ and the invertibility of $u_N$
modulo $\mu$ follow from $\gcd(N,\mu)=1$ and the identity:
\[
 u_N^2-N=(N^2+4)\left(\frac{N-1}{2}\right)^2.
\]
The normalization is fixed by
the odd quadratic Gauss sum:
\[
 \sum_{b\bmod\mu}\e(b^2/\mu)=\sqrt\mu,
 \qquad \mu\equiv1\pmod4.
\]
Both statements hold for composite $\mu$, including repeated
prime factors \cite[Section 30.1]{ER}.

For completeness, the full normalized matrix can be written
without any unspecified parameters. Relabel $X_{[a]}$ as
$D_{[b]}$ with $b=u_Na$, and put:
\[
 \alpha=\tfrac12\left(\frac1N+\frac1{\sqrt\mu}\right),\qquad
 \beta=\tfrac12\left(\frac1N-\frac1{\sqrt\mu}\right).
\]
In the order $(\one,T),(Y),(D)$, it is:
\begin{equation}
 S=\begin{pmatrix}
 \alpha&\beta&N^{-1}\mathbf j&-\mu^{-1/2}\mathbf j\\
 \beta&\alpha&N^{-1}\mathbf j& \mu^{-1/2}\mathbf j\\
 N^{-1}\mathbf j^{\mathsf t}&N^{-1}\mathbf j^{\mathsf t}&C_G&0\\
 -\mu^{-1/2}\mathbf j^{\mathsf t}&
 \mu^{-1/2}\mathbf j^{\mathsf t}&0&C_H
 \end{pmatrix},
 \label{eq:S}
\end{equation}
where each $\mathbf j$ is an all-ones row of the appropriate
length, and:
\[
 (C_G)_{[x,y],[x',y']}=\frac2N
       \cos\left(\frac{2\pi(xy'+x'y)}N\right),\qquad
 (C_H)_{[b],[b']}=\frac2{\sqrt\mu}
       \cos\left(\frac{4\pi bb'}\mu\right).
\]
The twists in the same order are:
\[
 T^{\mathrm{rib}}
 =\operatorname{diag}\bigl(1,1,\e(xy/N),\e(b^2/\mu)\bigr).
\]
In particular, the unit is fixed, and
$S_{\one,\one}=1/[N(1+t^2)]$.
Poisson summation on the refined frame grids gives both cosine
blocks and the vanishing mixed block. The spherical dimensions
give the first row. For the second row, let
$\mathcal E=(\DD_N)^0_A$ and let $b^0$ be its restriction
matrix. The local-unit row of $b^0$ is $\one+T$.
Taking this row in the normalized branching identity
$S^{\mathcal E}b^0=b^0S^{\DD_N}$ gives
$S_{\one,j}+S_{T,j}$ equal to $1/N$
for $j=\one,T$, to $2/N$ for $j=Y$, and to zero for $j=D$.
Thus the two distinguished rows follow directly from branching
\cite[Theorem 30.2]{ER}.
Under $X_{[a]}\mapsto D_{[u_Na]}$, the comparison with
\cite[Theorem 5.1 and Appendix F.1]{GSY} matches every entry of
$S$ and every ribbon twist, fixing the tensor unit.

When $\mu$ is prime, the spherical sign already determines the
metric group up to isometry
\cite[Corollary 5.6]{GSY}. At composite $\mu$, the sign alone
does not do so. The calculation here supplies the missing
parameter explicitly: for example, it selects the cyclic form
$b^2/85$ at $N=9$, and the cyclic form $b^2/125$ at $N=11$.
It does not assign that parameter to every independently
constructed HI category of the same order.

The comparison is between ribbon modular data. For ordinary
VOA characters the translation matrix includes the central-charge
scalar:
\[
 T^{\mathrm{an}}_{E_N}=\e(-c(E_N)/24)T^{\mathrm{rib}}
                    =\e(2N/3)T^{\mathrm{rib}}.
\]
For $W_N$, the pointed factor must also be retained.
Moreover, its analytic $S$-matrix is the negative of the
tensor-product categorical matrix with positive vacuum
normalization \cite[Sections 23.1 and 30]{ER}.

\section{Scope and related work}
\label{sec:scope}

There are three different identification questions. The
condensation proves a ribbon equivalence with the center of the
constructed $\HH_N$. The calculation in
Section~\ref{sec:modular} identifies its modular data with the
specified negative-branch GSY datum. An equivalence, or a
categorical Morita equivalence, with a separately prescribed
HI category requires further information.
A braided equivalence of centers would establish Morita
equivalence; a tensor equivalence of the underlying categories
requires stronger data, such as an equivalence of centers
preserving the canonical Lagrangian algebras.
Equality of modular data alone proves neither. The order-fifteen
examples of \cite[Corollary 5.9]{GSY} illustrate this failure
within the HI family.

A Galois automorphism with $\sqrt\mu\mapsto-\sqrt\mu$ takes
$-t$ to $\delta$ and yields a pseudo-unitary conjugate of
$\HH_N$ \cite[Theorem A(1) and Section 1.3]{GSY}.
This categorical operation supplies neither a VOA with that
conjugate representation category nor a positive-definite unitary
structure. A uniform VOA realization of the positive branch
therefore remains a separate problem.

Companion joint work \cite{Twisted} treats twisted
Haagerup--Izumi doubles on the negative spherical branch,
starting from a diagonal one-variable Virasoro family.
Its pointed associator is nontrivial; in the present cyclic
construction, the complementary Lagrangian factors give the
trivial associator. Both constructions use shifted $E_7$
lattice modules to cancel a pointed ribbon factor.
The present proof is independent of the forthcoming
twisted-case theorem.

\paragraph{Acknowledgments.}
The work of the first author was supported in part by a grant of
access to OpenAI models through the ChatGPT for Academic Researchers
programme. The second author is grateful to Pavel Etingof and Terry Gannon for their feedback and support.

\Needspace{12\baselineskip}

\end{document}